# Squarefree Integers And Extreme Values Of Some Arithmetic Functions
N. A. Carella, March, 2011.

**Abstract:** A study of the arithmetic function $\psi(N)/N = \prod_{p|N}(1+1/p)$, called the Dedekind psi function, concludes that its extreme values are supported on the subset of primorial integers $N = 2^{v_1} \cdot 3^{v_2} \cdots p_k^{v_k}$, where $p_i$ is the $k$th prime, and $v_i \geq 1$. In particular, the inequality $\psi(N)/N > 6\pi^{-2} e^\gamma \log\log N$ holds for all large integers $N = 2 \cdot 3 \cdot 5 \cdots p_k$, unconditionally. A discussion of the subsets of local minima and local maxima is also included.

## 1 Introduction

The psi function $\psi(N) = N\prod_{p|N}(1+1/p)$ and its normalized counterpart $\psi(N)/N = \prod_{p|N}(1+1/p)$ arise in various mathematic, and physic problems, see [SL]. Moreover, this function is entangled with other arithmetic functions. The values of the psi function coincide with the values of the sum of divisor function $\sigma(N) = \sum_{d|N} d$ on the subset of square-free integers. The product of the complex-valued psi function $\psi_s(N)/N^s = \prod_{p|N}(1+1/p^s)$ and the complex-valued totient function $\varphi_s(N)/N^s = \prod_{p|N}(1-1/p^s)$ satisfies the important relation

$$\frac{\psi_s(N)}{N^s} \frac{\varphi_s(N)}{N^s} = \prod_{p|N}(1-1/p^{2s}), \tag{1}$$

where $s \in \mathbb{C}$. The limit of this product as the integer $N_k = 2 \cdot 3 \cdot 5 \cdots p_k$ tends to infinity, where $p_k$ is the $k$th prime, is the zeta function $\zeta(2s) = \sum_{n \geq 1} n^{-2s}$. This establishes a clear-cut link to the Riemann hypothesis. This hypothesis claims that the nontrivial zeros of $\zeta(s)$ are located on the critical line $\{\ \mathfrak{Re}(s) = 1/2\ \}$, see [ES]. A survey of various strikingly different reformulations of the Riemann hypothesis appears in [AM].

The lower estimate given below appears to be new. A recent work [SP] claims that this estimate is equivalent to the Riemann hypothesis.

**Theorem 1.** Let $N \in \mathbb{N}$ be a primorial integer, then $\psi(N)/N > 6\pi^{-2} e^\gamma \log\log N$ holds unconditionally for all sufficiently large $N = 2 \cdot 3 \cdot 5 \cdots p_k$.

This lower estimate is quite similar to the lower estimate of the totient function, namely,

$$N_k / \varphi(N_k) > e^{\gamma} \log \log N_k \tag{2}$$

for all large integers $N_k$, see [NS], and [CL]. The reformulation of the Riemann hypothesis in terms of the psi function can be derived from (1) and (2), or by other means.

Two proofs of Theorem 1 are given in the last Section. The first proof, in Section 3.1, is obtained by an amalgam of elementary methods, but it is long. The second proof, in Section 3.2, relies on slightly more advanced results on the Prime Number Theorem, but it is significantly shorter. The next Section discussed background information.

## 2 Elementary Materials

The basic concepts and results employed throughout this work are stated in this Section.

**2.1. Density of Squarefree Integers.** Let $N = p_1^{v_1} \cdot p_2^{v_2} \cdots p_t^{v_t} \in \mathbb{N}$ be an integer, $v_i \geq 1$, $1 \leq i \leq t$. Let $\omega(n) = t$, and $\Omega(n) = v_1 + \text{L} + v_t$ be the prime divisors counting functions respectively. The Mobius function $\mu : \mathbb{R} \to \mathbb{Z}$ is defined by the relation

$$\mu(n) = \begin{cases} (-1)^{\omega(n)} & \text{if } \omega(n) = \Omega(n), \\ 0 & \text{if } \omega(n) \neq \Omega(n). \end{cases} \tag{3}$$

This function is identically zero on the set of nonsquarefree integers, so $\mu^2(n) = |\mu(n)|$ is the characteristic function of squarefree numbers. The Dirichlet average of the Mobius function over the divisor of an integer has the simple expression

$$\sum_{d|n} \mu(n) = \begin{cases} 1 & \text{if } n = 1, \\ 0 & \text{if } n \neq 1. \end{cases} \tag{4}$$

**Theorem 2.** The number of squarefree integers in the interval $[1, x]$ is $Q(x) = \sum_{n \leq x} |\mu(n)| = 6\pi^{-2}x + O(x^{1/2})$.

Proof: Every integer $n \in \mathbb{N}$ has a unique representation as $n = ab^2$, where $a$ is a squarefree integer. This means that the characteristic function of the squarefree integer can be rewritten as

$$\mu^2(n) = \sum_{d^2 | n} \mu(d), \tag{5}$$

and the number of squarefree integers up to $x$ is $\sum_{n \leq x} |\mu(n)| = \sum_{n \leq x} \mu(n)^2 = \sum_{n \leq x} \sum_{d^2 | n} \mu(d)$. So the total number of squarefree integers $n \leq x$ is given by

$$\sum_{n \leq x} \mu^2(n) = \sum_{n \leq x} \sum_{d^2 | n} \mu(d) = \sum_{d < x^{1/2}} \left[ \frac{x}{d^2} \right] \mu(d), \tag{6}$$

where the bracket $[\ x\ ]$ denotes the largest integer function, and $[x/d^2]$ tallies the number of integers $n \leq x$, and divisible by $d^2$. Now remove the bracket in the previous expression to obtain

$$\sum_{d<x^{1/2}} \left[\frac{x}{d^2}\right] \mu(d) = \frac{6x}{\pi^2} + O\left(x \sum_{d \geq x^{1/2}} d^{-2} + \sqrt{x}\right). \tag{7}$$

The last step follows from the identity $6/\pi^2 = \sum_{n \geq 1} \mu(n) n^{-2} = \sum_{n \leq x^{1/2}} \mu(n) n^{-2} + \sum_{n > x^{1/2}} \mu(n) n^{-2})$. ∎

Since it is known that $\sum_{n \leq x} |\mu(n)| = 6\pi^{-2} x + \Omega_{\pm}(x^{1/4})$, see [TN, p. 128], the best error term possible seems to be at least $O(x^{1/4+\varepsilon})$. Trivially, the short interval $(x, x + y]$ contains $\sum_{x < n \leq x+y} |\mu(n)| = 6\pi^{-2} y + O(x^{1/2})$ squarefree integers for all $y \geq cx^{1/2}$, $c > 0$ constant. A better and sharper estimate is available in the literature.

An interesting application of the squarefree counting measure $d(6\pi^{-2} x + O(x^{1/2}))$ is in estimating the squarefree harmonic sum.

***Theorem 3.*** Let $x \geq x_0$ be a large real number. Then $\sum_{n \leq x} \mu^2(n) n^{-1} = 6\pi^{-2} \log x + O(x^{-1/2})$.

Proof: Put $f(x) = 1/x$. Then, integral representation of the harmonic sum over the squarefree integers is given by

$$\sum_{n \leq x} \frac{\mu^2(n)}{n} = \int_c^x \frac{d(6\pi^{-2} t + O(t^{1/2}))}{t}, \tag{8}$$

where $c = e^{6/\pi^2}$. Integrating by parts yields the assertion. ∎

Similarly, exponential sums, and other finite sums over the squarefree integer can be approximated via the integral representations. For example,

$$\sum_{n \leq x} \mu^2(n) \chi(n) e^{i2\pi a n} = \int_c^x \chi(t) e^{i2\pi a t} d(6 e^\gamma \pi^{-2} t + O(t^{-1/2})), \tag{9}$$

where $\chi$ is a character modulo $q \geq 2$.

**2.2 Sums and Product Over the Primes.** The most basic finite sum over the prime numbers is the prime harmonic sum $\sum_{n \leq x} p^{-1}$. The standard evaluation of this sum is stated below.

***Theorem 4.*** (Mertens) Let $x \geq 2$. Then

(i) $\sum_{p \leq x} \frac{1}{p} = \log \log x + B_1 + O(1/\log^2 x)$, (ii) $\sum_{p \leq x, p \equiv a \bmod q} \frac{1}{p} = \frac{1}{\varphi(q)} \log \log x + B_{a,q} + O_q(1/\log^2 x)$,

where $B_1 = .2614972128\ldots$, and $B_{a, q}$ is Mertens constant for the arithmetic progression $qn + a$.

The proofs of these results, and the simpler and elementary estimate

$$\left|\sum_{p \leq x} 1/p - \log \log x + B_1\right| \leq \frac{c}{\log x} \tag{10}$$

where $c > 0$ is a constant are widely available in the literature, see [TN]. There are various versions of the error term $R(x)$ in $\sum_{n \leq x} p^{-1} = \log \log x + B_1 + R(x)$. Currently, the best explicit estimates of the error term are

$$R(x) = \pm\left(\frac{1}{10\log^2 x} + \frac{4}{15\log^3 x}\right) \quad \text{or} \quad R(x) = \pm\frac{3\log x + 4}{8\pi\sqrt{x}} \tag{11}$$

unconditionally, see [DU], or modulo the Riemann hypothesis respectively, see [RS], [VL. p. 26]. There are various analytical formulas for the constant $B_1$, and Euler constant $\gamma = \lim_{x\to\infty} \sum_{n\leq x}(n^{-1} - \log n)$, confer the literature. The one given below will be used later on.

***Lemma 5.*** The constant $B_1 = \gamma - \sum_{p\geq 2}\sum_{n\geq 2}(np^n)^{-1}$.

Proof: The Mertens constant is defined by $B_1 = \gamma + \sum_{p\geq 2}(\log(1-1/p) + 1/p)$, see [HW, p. 466]. Using the power series $\log(1+x) = \sum_{n\geq 1}(-1)^{n+1}x^n/n$, $0 \leq x < 1$, and simplifying, return the required expression. ∎

As an application of the penultimate result, there are the following interesting products.

***Theorem 6.*** Let $x \in \mathbb{R}$ be a large real number, then

$$\prod_{p\leq x}(1-1/p)^{-1} = e^\gamma \log x + O(1/\log^2 x). \tag{12}$$

Proof: This follows from Theorem 4. ∎

***Theorem 7.*** Let $x \geq x_0$ be a real number, then

$$\prod_{p\leq x}(1+1/p) = 6\pi^{-2}e^\gamma \log x + O(1/\log^2 x). \tag{13}$$

Proof: For a large real number $x \in \mathbb{R}$, rewrite the product as

$$\prod_{p\leq x}(1+1/p)(1-1/p)(1-1/p)^{-1} = \prod_{p\leq x}(1-1/p^2)\prod_{p\leq x}(1-1/p)^{-1}. \tag{14}$$

Replacing $\prod_{p\geq 2}(1-1/p^2) = \sum_{n\geq 1}\mu(n)n^{-2} = 6\pi^{-2}$ in the first product on the right side, yields

$$\prod_{p\leq x}(1-1/p^2)\prod_{p\leq x}(1-1/p)^{-1} = \left(\frac{6}{\pi^2} - \sum_{n\geq z}\mu(n)n^{-2}\right)\prod_{p\leq x}(1-1/p)^{-1}$$
$$= \frac{6}{\pi^2}\prod_{p\leq x}(1-1/p)^{-1} + O(\frac{\log x}{x^c}), \tag{15}$$

where $z = O(x^c)$, $c > 1$ constant. Lastly, applying Theorem 6, to the last product above, yields the claim. ∎

This result immediately gives the maximal value of the sum of divisors function at the squarefree numbers, viz,

$$\frac{\sigma(N)}{N} \leq \frac{6e^\gamma}{\pi^2}\log\log N + O(\frac{1}{(\log\log N)^2}) \tag{16}$$

for any squarefree integer $N \geq 1$.

**Theorem 8.** ([DP])  The quantity $x^{1/2}\left(\prod_{p\leq x}(1-1/p)^{-1} - e^{\gamma}\log x\right)$ attains arbitrary large positive and negative values as $x \to \infty$.

**2.3. Squarefree Integers And Prime Numbers Link.** The asymptotic density of squarefree integers in the interval $[1, x]$ is denser than the density of prime numbers by a factor of $(\pi^2/6)\log x$. Nevertheless, some results on sums over squarefree integers are closely linked to products over the primes. Some of these results are asymptotically related by a factor of $e^{\pm\gamma}$. The factor $e^{\pm\gamma}$ has been a topic of contention for quite sometime, see [WG, p. 107] for a discussion of its link to the Sieve of Eratosthenes, [MI] for its link to the density of primes in short intervals, [WF] for its link to Gillies conjecture, et cetera. For example, asymptotically, the squarefree harmonic sum and the prime harmonic product are equivalent up to the constant factor $e^{\pm\gamma}$. A comparison of these two expressions, shown below side by side, reveals the strikingly different error terms:

$$\sum_{n\leq x}\mu^2(n)n^{-1} = 6\pi^{-2}\log x + O(x^{-1/2}) \quad \text{and} \quad \prod_{p\leq x}(1+1/p) = 6\pi^{-2}e^{\gamma}\log x + O(1/\log^2 x). \tag{17}$$

Let $P_x = \prod_{p\leq x} p$, and let $\kappa = e^{\gamma}$. It is quite easy to see that

$$\sum_{n\leq x}\frac{\mu^2(n)}{n} = \sum_{d\mid P_x}\frac{1}{d} - \sum_{x<d\mid P_x}\frac{1}{d}, \tag{18}$$

where $\sum_{d\mid P_x}1/d = \prod_{p\leq x}(1+1/p)$. But it seems to be difficult to show that

$$\sum_{n\leq x^{\kappa}}\frac{\mu^2(n)}{n} + \sum_{x<d\mid P_x}\frac{1}{d} \geq \sum_{d\mid P_x}\frac{1}{d}. \tag{19}$$

Specifically, confirming the upper estimate

$$\sum_{x<d\mid P_x}\frac{1}{d} \leq c\log x + O(\log x/\sqrt{x}) \tag{20}$$

is a difficult problem. An affirmative answer will immediately leads to the optimum error term bound in the product $\prod_{p\leq x}(1+1/p) = 6\pi^{-2}e^{\gamma}\log x + O(\log x/\sqrt{x})$, as specified by the Riemann hypothesis, see [RS].

**3 Squarefree Integers And The Extreme Values of Arithmetic Functions**
The first proof of Theorem 1, presented below, is completely elementary, but longer. In comparison, the second proof relies on slightly more advanced methods, but it is significantly shorter.

Let $N \in \mathbb{N}$ be an integer. The sum of divisors function $\sigma(N)/N = \sum_{d\mid N}1/d$ restricted to the subset of squarefree integers is called the Dedekind psi function and denoted by $\psi(N)/N = \prod_{p\mid N}(1+1/p)$. The psi function is an oscillatory function on the interval $[1, x]$, its values fall in the range

$$1 + \frac{1}{p_{\max}} \leq \frac{\psi(N)}{N} \leq \left(1+\frac{1}{2}\right)\left(1+\frac{1}{3}\right)\cdots\left(1+\frac{1}{p_k}\right), \tag{21}$$

where $p_{\max} \leq x$, and $2 \cdot 3 \cdots p_k \leq x$ are the largest prime, and the largest primorial integers in the interval $[1, x]$ respectively. The restriction $\psi(N_k)/N_k$ to primorial integers $N_k$ is strictly monotonically increasing as $N_k \to \infty$. The jump-discontinuity $\Delta_k$ at the point $N_{k+1}$ has the magnitude

$$\Delta_k = \frac{\psi(N_{k+1})}{N_{k+1}} - \frac{\psi(N_k)}{N_k} = \frac{\psi(N_k)}{p_{k+1}}, \qquad (22)$$

where $\dfrac{1}{p_k + p_k^\alpha} < \dfrac{1}{p_{k+1}} \leq \dfrac{1}{p_k + 2}$, $0 < \alpha < .526$, see [BK]. The quantity $\Delta_k$ encodes important information on the distribution of prime numbers. For example, a pair of twin primes $p_k$ and $p_{k+1}$ produces the largest jump-discontinuity $\Delta_k = \psi(N_k)/(p_k + 2)$ possible.

**3.1 First Proof.** From the inequality

$$\prod_{p \mid N}(1 + 1/p) = \prod_{p \mid N}(1 - 1/p^2) \prod_{p \mid N}(1 - 1/p)^{-1} > \frac{6}{\pi^2} \prod_{p \mid N}(1 - 1/p)^{-1}, \qquad (23)$$

and Theorem 6, it is plausible to expect that $\prod_{p \mid N}(1 + 1/p) > 6e^\gamma \pi^{-2} \log \log N$ infinitely often as $N \to \infty$. This fact is proven below by elementary methods.

***Theorem 1*** Let $p_k$ be the $k$th prime, and let $N_k = 2 \cdot 3 \cdot 5 \cdots p_k$, $k \geq 1$. Then $\psi(N_k)/N_k > 6\pi^{-2} e^\gamma \log \log N_k$ for all sufficiently large integer $N_k$.

Proof: Assume the claim is false, that is, $\psi(N_k)/N_k = \prod_{p \mid N_k}(1 + 1/p) \leq 6\pi^{-2} e^\gamma \log \log N_k$ for $N_k \geq N_0$. Then

$$\log \prod_{p \mid N_k}(1 + 1/p) \leq \log\left(6\pi^{-2} e^\gamma \log \log N_k\right). \qquad (24)$$

Expanding the left side into power series and separating into linear and nonlinear sums yield

$$\sum_{p \mid N_k} \log(1 + 1/p) = \sum_{p \leq p_k} \sum_{n=1}^{\infty} \frac{(-1)^{n+1}}{np^n} = \sum_{p \leq p_k} \sum_{n=2}^{\infty} \frac{(-1)^{n+1}}{np^n} + \sum_{p \leq p_k} \frac{1}{p}, \qquad (25)$$

since $p_k \mid N_k$ is equivalent to $p \leq p_k$. Applying Theorem 4, to the previous expression transforms it into

$$\sum_{p \leq p_k} \log(1 + 1/p) = \sum_{p \leq p_k} \sum_{n=2}^{\infty} \frac{(-1)^{n+1}}{np^n} + \log \log p_k + B_1 + O\left(\frac{1}{(\log \log p_k)^2}\right), \qquad (26)$$

where $B_1$ is a constant. Next replace $B_1 = \gamma - \sum_{p \geq 2} \sum_{n \geq 2} (np^n)^{-1}$, see Lemma 5, in the previous equation to obtain

$$\sum_{p \mid N_k} \log(1 + 1/p) = \sum_{p \leq p_k} \sum_{n=2}^{\infty} \frac{(-1)^{n+1}}{np^n} - \sum_{p \geq 2,\ n=2}^{\infty} \frac{1}{np^n} + \log \log p_k + \gamma + O\left(\frac{1}{(\log \log N_k)^2}\right). \qquad (27)$$

Now, inequality (24) becomes

$$\sum_{p \leq p_k} \sum_{n=2}^{\infty} \frac{(-1)^{n+1}}{np^n} - \sum_{p \geq 2,} \sum_{n=2}^{\infty} \frac{1}{np^n} + \log \log p_k + \gamma + O(\frac{1}{(\log \log N_k)^2}) \leq \gamma + \log(6/\pi^2) + \log \log \log N_k. \tag{28}$$

This is equivalent to the inequality

$$\log \log p_k - \log \log \log N_k + O(\frac{1}{(\log \log N_k)^2}) \leq \log(6/\pi^2) + \sum_{p \leq p_k} \sum_{n=2}^{\infty} \frac{(-1)^n}{np^n} + \sum_{p \geq 2,} \sum_{n=2}^{\infty} \frac{1}{np^n}. \tag{29}$$

A routine exercise using $\vartheta(p_k) = p_k + o(p_k) \leq c_0 \log N_k$, $\log N_k = \sum_{p \leq p_k} \log p = \vartheta(p_k)$, and $p_k \leq c_0 \log N_k$, $c_0 > 1$ constant, see [CL, Proposition 6], demonstrates that the logarithmic difference

$$\log \log p_k - \log \log \log N_k = O(\frac{1}{(\log \log N_k)^{B+1}}), \tag{30}$$

where $B > 0$ is an arbitrary real number. Combining (29) and (30) yield

$$O(\frac{1}{(\log \log N_k)^{B+1}}) + O(\frac{1}{(\log \log N_k)^2}) \leq \log(6/\pi^2) + \sum_{p \leq p_k} \sum_{n=2}^{\infty} \frac{(-1)^n}{np^n} + \sum_{p \geq 2,} \sum_{n=2}^{\infty} \frac{1}{np^n}. \tag{31}$$

The right side of the previous inequality can be simplified as

$$\log(6/\pi^2) + \sum_{p \leq p_k} \sum_{n=2}^{\infty} \frac{(-1)^n}{np^n} + \sum_{p \geq 2, n=2}^{\infty} \frac{1}{np^n} = -\sum_{p \geq 2} \sum_{n=1}^{\infty} \frac{1}{np^{2n}} + \sum_{p \leq p_k} \sum_{n=2}^{\infty} \frac{(-1)^n}{np^n} + \sum_{p \geq 2, n=2}^{\infty} \frac{1}{np^n}$$

$$= -\sum_{p > p_k} \sum_{n=1}^{\infty} \frac{1}{2np^{2n}} + \sum_{p > p_k} \sum_{n=1}^{\infty} \frac{1}{(2n+1)p^{2n+1}} \tag{32}$$

$$< 0.$$

This last regrouping stems from the expansion $\log(6/\pi^2) = \log \prod_{p \geq 2}(1 - 1/p^2) = -\sum_{p \geq 2} \sum_{n \geq 1}(np^{2n})^{-1}$. Replacing the previous one into the penultimate one returns

$$O(\frac{1}{(\log \log N_k)^{B+1}}) + O(\frac{1}{(\log \log N_k)^2}) < 0, \tag{33}$$

where $B > 1$, see [SH, p. 451]. Now, observe that for any primorial integer $N_k = 2 \cdot 3 \cdot 5 \cdots p_k$, the error term of the prime harmonic sum

$$\sum_{p \leq p_k} \frac{1}{p} = \log \log p_k + B_1 + O(\frac{1}{(\log \log p_k)^2}), \tag{34}$$

where $p_k \leq c_0 \log N_k$, $c_0 > 1$ constant, is nonnegative, so the left side of (33) is a nonnegative quantity that approaches 0 from the right. But the right side of (33) is a negative quantity that approaches 0 from the left as $N_k \to \infty$. Clearly, this is a contradiction for all sufficiently large integers $N_k > 1$ as claimed. ∎

This analysis also agrees with the well known asymptotic limit of the product of the reciprocal of the primes, viz,

$$\lim_{n \to \infty} \frac{1}{\log n} \prod_{p \mid n} (1 + 1/p) = \frac{6e^\gamma}{\pi^2}. \tag{35}$$

**3.2 Second Proof.** Another more advanced and compact approach to the proof of the inequality $\psi(N_k)/N_k > 6\pi^{-2} e^\gamma \log \log N_k$ for all large integer $N_k$, is investigated in this Section. The derivation of this result unfolds from the recent result on the oscillation theorem of finite prime product in [DP].

***Theorem* 1.** Let $N \in \mathbb{N}$ be a primorial integer, then $\psi(N)/N > 6\pi^{-2} e^\gamma \log \log N$ holds unconditionally for all sufficiently large $N = 2 \cdot 3 \cdot 5 \cdots p_k$.

Proof : Theorem 8 implies that the product

$$\prod_{p \leq x} (1 - 1/p)^{-1} = e^\gamma \log x + \Omega_\pm (\log \log x / \sqrt{x}), \tag{36}$$

where $f(x) = g(x) + \Omega_\pm(h(x))$ means that either $f(x) > g(x) + ch(x)$ or $f(x) < g(x) - ch(x)$ infinitely often, $c > 0$ constant. From the identity

$$\prod_{p \leq x}(1 + 1/p) = \prod_{p \leq x}(1 - 1/p^2) \prod_{p \leq x}(1 - 1/p)^{-1} > \frac{6}{\pi^2} \prod_{p \leq x} (1 - 1/p)^{-1}, \tag{37}$$

and (36), it follows that

$$\prod_{p \leq x}(1 + 1/p) > \frac{6e^\gamma}{\pi^2} \log x + \Omega(\log \log x / \sqrt{x}) \quad \text{and} \quad \prod_{p \leq x}(1 + 1/p) < \frac{6e^\gamma}{\pi^2} \log x - \Omega(\log \log x / \sqrt{x}). \tag{38}$$

It shows that $\prod_{p \leq x}(1 + 1/p)$ oscillates infinitely often symmetrically about the line $6e^\gamma \pi^{-2} \log x$ as $x \to \infty$. To rewrite the variable $x \geq 1$ in terms of the integer $N$, recall that

$$\log N_k = \sum_{p \leq p_k} \log p = \vartheta(p_k), \quad \text{and} \quad \vartheta(p_k) = p_k + o(p_k) \leq c \log N_k. \tag{39}$$

So it readily follows that $p_k \leq x = c \log N_k$, where $c > 1$ is a constant, and $\vartheta(x)$ is the Chebychev function. Moreover, since the maxima of the sum of divisor function $\sigma(N) \geq \psi(N)$ occur at the colossally abundant integers $N = 2^{v_1} \cdot 3^{v_2} \cdots p_k^{v_k}$, and $v_1 \geq v_2 \geq \cdots \geq v_k \geq 1$, see [AE], [BR], [LA], [RJ], it follows that the maxima of the Dedekind psi function $\psi(N)$ occur at primorial integer $N_k = 2 \cdot 3 \cdot 5 \cdots p_k$. Therefore, the expressions (38) and (39) imply that

$$\prod_{p \leq c \log N_k}(1 + 1/p) > \frac{6e^\gamma}{\pi^2} \log \log N_k + \Omega_+(\log \log \log N_k / \sqrt{\log N_k}). \tag{40}$$

as the primorial integer $N_k = 2 \cdot 3 \cdot 5 \cdots p_k$ tends to infinity. ∎

**3.3 The Local Minima and Maxima.** Define the subsets of integers

$$m_\psi = \{ n : \psi(n)/n < 6e^\gamma \pi^{-2} \log\log n \} \quad \text{and} \quad M_\psi = \{ n : \psi(n)/n > 6e^\gamma \pi^{-2} \log\log n \}. \tag{41}$$

The densities of the integers in these subsets are asymptotically the same as the density of prime numbers. More precisely,

$$M_\psi(x) = \#\{ n \leq x : \psi(n)/n > 6e^\gamma \pi^{-2} \log\log n \} = \frac{x}{\log x}\left(1 + O\left(\frac{1}{(\log\log x)^2}\right)\right), \tag{42}$$

the same estimate holds for $m_\psi(x)$. These densities results are derived from the distribution function

$$D_\alpha(t) = \exp\left(-e^{e^{-\gamma}\zeta(2)t}\left(1 + O(t^{-2})\right)\right) = \frac{\#\{n \leq x : \psi(n) > t\}}{x}, \tag{43}$$

where $x \geq x_0$, of the psi function, see [WN].